\documentclass[5pt]{article}
\usepackage{amsmath}
\usepackage[latin1]{inputenc}
\usepackage{amsmath,amsthm,amssymb}
\usepackage{graphicx}
\usepackage{color}

\renewcommand{\theequation}{\thesection.\arabic{equation}}

\newtheorem{lemma}{Lemma}[section]
\newtheorem{thm}{Theorem} [section]

\newtheorem{exmp}{Example} [section]
\newtheorem{cor}{Corollary}[section]
\newtheorem{rem}{Remark}[section]
\newtheorem{conj}{Conjecture}[section]
\title{On some conjectures of exponential Diophantine equations
\thanks {The work is supported by the Natural Science Foundation of China with No.11671153, No.11971180.}}
\author{{Hairong Bai \footnote{\ H.R.Bai, Institute of Mathematics Science, South China Normal University,
 \small Guangzhou, Guangdong, China (Email: Baihairong2007@163.com)}
}\ }
\date{}
\begin{document}
\baselineskip15pt \maketitle
\renewcommand{\theequation}{\arabic{section}.\arabic{equation}}
\catcode`@=11 \@addtoreset{equation}{section} \catcode`@=12
\begin{abstract}

In this paper, we consider the exponential Diophantine equation $a^{x}+b^{y}=c^{z},$ where $a, b, c$ be relatively prime positive integers such that $a^{2}+b^{2}=c^{r}, r\in Z^{+}, 2\mid r$ with $b$ even. That is
$$a=\mid Re(m+n\sqrt{-1})^{r}\mid, b=\mid Im(m+n\sqrt{-1})^{r}\mid, c=m^{2}+n^{2},$$ where $m, n$ are positive integers with $m>n, m-n\equiv1(mod 2),$ gcd$(m, n)=1.$ $(x, y, z)= (2, 2, r)$ is called the trivial solution of the equation. In this paper we prove that the equation has no nontrivial solutions in positive integers $x, y, z$ when
$$r\equiv 2(mod 4), m\equiv 3(mod 4), m>\max\{n^{10.4\times10^{11}\log(5.2\times10^{11}\log n)}, 3e^{r}, 70.2nr\}.$$
Especially the equation has no nontrivial solutions in positive integers $x, y, z$ when
$$r=2, m\equiv 3(mod 4), m>n^{10.4\times10^{11}\log(5.2\times10^{11}\log n)}.$$


 \end{abstract}

 \section{Introduction and the main result}
 $ $

In this paper we consider the exponential Diophantine equation $$a^{x}+b^{y}=c^{z} \eqno{(1.1)},$$ where $a, b, c$ be relatively prime positive integers.

Equation (1.1) was first considered by Mahler ~\cite{Mah33}. For a fixed triple $(a,b,c),$ he proved that (1.1) has finitely many solutions $(x, y, z)$ under the assumption that $a, b, c>1.$ Since his method was based on a p-adic generalisation of the Thue-Siegel method, it was ineffective. Later, Gel$^{,}$fond ~\cite{Ge40} gave an effective result; his method was based on Baker$^{\prime}$s theory of linear forms in logarithms of algebraic numbers.

Using elementary number theory methods such as congruences, Jacobi symbols and standard divisibility arguments involving ideals in quadratic (or cubic) number fields, all solutions of (1.1) where a,b,c are distinct primes $\leq17$ were determined (see ~\cite{DYL18}, ~\cite{Ha76}, ~\cite{Na58} and ~\cite{Uc76}).

In 1956, Je\'{s}smanowicz conjectured that:
\begin{conj}
If $a^{2}+b^{2}=c^{2},$ then the exponential Diophantine equation $a^{x}+b^{y}=c^{z}$ has no nontrivial solutions in positive integers $x, y, z$ except for the trivial solution $(x, y, z)= (2, 2, 2).$
\end{conj}

If $a^{2}+b^{2}=c^{2},$ then $$a=m^{2}-n^{2}, b=2mn, c=m^{2}+n^{2},$$ furthermore $$\max\{a, b\}<c<\min\{a^{2}, b^{2}\}.$$

Sierpi\'{n}ski~\cite{Si55} proved that if $(m, n)=(2, 1),$ then conjecture 1.1 is true. Later, Je\'{s}manowicz~\cite{Je55} proved that if $(m, n)=(3, 2), (4, 3), (5, 4), (6, 5),$ then conjecture 1 is true.

Many author has proved that conjecture 1.1 is true for Various concrete values of $n$ without any restrictions on $m.$ Lu~\cite{Lu59} $n=1.$  Demjanenko~\cite{De65} proved that if $m=n+1,$ then conjecture 1.1 is true.

N. Terai~\cite{Te14} $n=2,$ T. Miyazaki and N. Terai~\cite{MiT15} $n<100, ord_{2}(n)=1.$

An extension of the Je\'{s}manowicz conjecture, was suggested by Terai. He conjectured that
\begin{conj}
If $a^{2}+b^{2}=c^{r}, r\geq 2,$ then the exponential Diophantine equation $a^{x}+b^{y}=c^{z}$ has no nontrivial solutions in positive integers $x, y, z$ except for the trivial solution $(x, y, z)= (2, 2, r).$
\end{conj}

Let $a, b, c$ be relatively prime positive integers such that $a^{2}+b^{2}=c^{r}, r\in Z^{+}, 2\mid r$ with $b$ even. That is
$$a=\mid Re(m+n\sqrt{-1})^{r}\mid, b=\mid Im(m+n\sqrt{-1})^{r}\mid, c=m^{2}+n^{2},$$ where $m, n$ are positive integers with $m>n, m-n\equiv1(mod 2),$ gcd$(m, n)=1.$

For the case $2\nmid n,$ T. Miyazaki~\cite{Mi18} proved that, Suppose $r\equiv 2(mod 4)$ and
\  \begin{displaymath}\begin{cases}
    \ n\equiv 3(mod 4),  &  \text{$r\equiv 2(mod 8)$}, \\
    \ n\equiv 1(mod 4),  &  \text{$r\equiv 6(mod 8)$}, \\
 \end{cases}                \end{displaymath}
Then there exists an effectively computable constant $C,$ which depends only on $r$ and $n,$ such that Conjecture 1.2 is true whenever $m>C.$
Especially, for any fixed $n$ with $n\equiv3(mod 4),$ conjecture 1.1 is true except for finitely many values of $m.$

In this paper we consider the case $2\mid n,$ then we may assume that $n\geq 4.$ We establish the following results.

\begin{thm}\label{thm1}
Suppose $r\equiv 2(mod 4), m\equiv 3(mod 4).$ If $m>\max\{n^{10.4\times10^{11}\log(5.2\times10^{11}\log n)}, 3e^{r}, 70.2nr\},$ then Conjecture 1.2 is true.
\end{thm}

\begin{cor}\label{thm2}
Suppose that $m\equiv 3(mod 4).$ If $m>n^{10.4\times10^{11}\log(10^{11}\log n)},$ then Conjecture 1.1 is true.
\end{cor}

The result [Mi18, Theorem 1.1] may be refined:
\begin{rem}\label{thm1}
Suppose $r\equiv 2(mod 4)$ and
\  \begin{displaymath}\begin{cases}
    \ n\equiv 3(mod 4),  &  \text{$r\equiv 2(mod 8)$}, \\
    \ n\equiv 1(mod 4),  &  \text{$r\equiv 6(mod 8)$}, \\
 \end{cases}                \end{displaymath}
If $m>\max\{n^{10.4\times10^{11}\log(5.2\times10^{11}\log n)}, 3e^{r}, 70.2nr\},$ then Conjecture 1.2 is true. Especially if $n\equiv 3(mod 4), m>n^{10.4\times10^{11}\log(10^{11}\log n)},$ then Conjecture 1.1 is true.
\end{rem}

Let $n=2^{\alpha}i, m=2^{\beta}j+e, e\in\{1, -1\}, 2\nmid ij.$ Then $\alpha\geq1, \beta\geq 2.$

Let $n_{(p)}=p^{ord_{p}(n)},$ where $ord_{p}(n)=\max\{\tau, p^{\tau}\mid n\}.$

\begin{thm}\label{thm1}
Suppose that $r=2, m\equiv 3(mod 4), m>56n.$ If one of the follow cases holds, then Conjecture 1.1 is true.

(1)\ $(m, n)$ are not the case $2\alpha=\beta+1, j\equiv 1(mod 4), \alpha\geq 3.$ Furthermore $$n^{\nu}<4^{\alpha+1}, m>\sqrt{n^{\eta}-n^{2}},$$ where $$\nu=1+\frac{\log 2}{\log mn}, \eta=\frac{(22\nu+1)(\log 2)}{11\log 4^{\alpha+1}/n^{\nu}}.$$

(2)\ $n=2^{\alpha}.$
\end{thm}

\begin{thm}\label{thm1}
Suppose $r=2, m\equiv 3(mod 4).$ Assume that $m>56n$ or $$\max\{n_{(3)}, n_{(5)}\}\geq (1534\log n)^{1.5}t(n)\sqrt{n}.$$
If one of the follow case holds, then Conjecture1.1 is true.

(1)$(m, n)$ are not the case $2\alpha=\beta+1, j\equiv 1(mod 4), \alpha\geq 3.$ (i) holds.

(2)$(m, n)$ are the case $2\alpha=\beta+1, j\equiv 1(mod 4).$ (i), (ii) hold.

Where (i) $2n_{(2)}n_{(3)}n_{(5)}\geq 2^{s(n)}\sqrt{n}.$

(ii) $\max\{n_{(3)}, n_{(5)}\}\geq 2\sqrt{n\log n}$ and $n_{(3)}n_{(5)}\geq\frac{3^{1/12}}{2^{\alpha+1}}n.$
$$s(n)=\frac{2\log(n^{11}-4n^{2})+5\log n+2\log(n^{11}/4)}{2(\log(n^{11}/4))\log((n^{11}-4n^{2})n^{2}/4)-2(\log2)\log n}\log n(<1.21),$$
$$t(n)=\sqrt{\frac{1}{1534}+(1+\frac{4}{n})\log(1+\frac{2n^{2}}{6n+9})}\frac{\log(1+\frac{2n^{2}}{6n+9})}{(\log n)^{1.5}}(<\sqrt{2}).$$
\end{thm}

\begin{cor}\label{thm1}
Suppose $r=2, m\equiv 3(mod 4)$ and that $(m, n)$ not the case $2\alpha=\beta+1, j\equiv 1(mod 4), \alpha\geq 3.$ If one of the follow cases holds, then Conjecture 1.1 is true. (1)\ $n\leq 64.$ (2)\ $n>64, n_{(2)}>\sqrt{n}, m>56n.$
\end{cor}

\begin{cor}\label{thm1}
If $r=2, m\equiv 3(mod 4)$ and $n\in\{6^{u}, u\geq94,  10^{v}, v\geq 40\}, u, v\in Z^{+},$ then Conjecture 1.1 is true.
\end{cor}

In this paper, we always assume that $(x, y, z)$ is a solution of the equation (1.1).
we often consider the case $$r=2 \quad or \quad r>2, r\equiv 2(mod 4), m>\frac{2r}{\pi}n \eqno{(1.2)}$$

In section 2, we will give some congruence properties of solutions of (1.1). Let $$\Delta=\mid \frac{r}{2}x-z\mid, K=\frac{m}{n}, F_{r}(K)=\frac{c^{r/2}}{a}.$$ In section 3, we will give estimates of $Y, \Delta$ and $K.$ In section 4, we consider the case $y=1.$
In section 5, we prove Theorem 1.1 and Corollary 1.1.
In section 6, we prove Theorem 1.2 and Theorem 1.3.
In section 7, we prove Corollary 1.2. and Corollary 1.3.
In section 8, we will give estimates of $s(n), t(n)$ and an example of the case of $m=\prod_{i=1}^{k}p_{i}^{l_{i}}, p_{i}\equiv 1(mod 4).$

\section{The congruence properties of solutions}
$ $

We always assume that $(x, y, z)$ is a solution of the equation (1.1). In this section, we will give some congruence properties of solutions.

Write $m+n\sqrt{-1}=c^{1/2}e^{i\psi}$ with $\psi\in(0, \pi/2).$ Then $r\psi=r\arctan(n/m)<rn/m.$ If the case (1.2) holds, Then $r\psi<rn/m<\pi/2.$
Hence $Re(m+n\sqrt{-1})^{r}, Im(m+n\sqrt{-1})^{r}$ are positive. So $$a=m^{r}-C_{r}^{2}m^{r-2}n^{2}+\cdots+C_{r}^{r-2}m^{2}n^{r-2}-n^{r},$$
$$b=rm^{r-1}n-C_{r}^{3}m^{r-3}n^{3}+\cdots+C_{r}^{r-3}m^{3}n^{r-3}-rmn^{r-1}.$$

Especially If $a^{2}+b^{2}=c^{2},$ Then $$a=m^{2}-n^{2}, b=2mn, c=m^{2}+n^{2},\quad \max\{a, b\}<c<\min\{a^{2}, b^{2}\}.$$
\begin{lemma}\label{lemma1}
Suppose that (1.2) holds.

(1)\ If gcd$(m, n^{2}-1)>1,$ then $2\mid x.$

(2)\ If gcd$(m, n^{2}+1)\notin\{1, 2^{t}\},$ then $2\mid z.$

(3)\ If gcd$(m^{2}+1, n)\notin\{1, 2^{t}\},$ then $x\equiv z(mod 2).$
\end{lemma}

Proof. (1)\ If $2\mid m,$ then $2\nmid n.$ Hence
$$a\equiv -n^{r}\equiv -1(mod4), b\equiv rmn^{r-1}\equiv 0(mod4), c=m^{2}+n^{2}\equiv 1(mod4).$$ It is easy to see that $2|x$ by taking (1.1) modulo 4.

If $2\nmid m,$ gcd$(m, n^{2}-1)>1.$ Then $2\nmid gcd(m, n^{2}-1).$ Suppose that $p>2$ is a prime and $p\mid gcd(m, n^{2}-1).$
Hence $$a\equiv -n^{r}\equiv -(n^{2})^{r/2}\equiv -1(mod p), b\equiv 0(mod p), c=m^{2}+n^{2}\equiv 1(mod p).$$ It is easy to see that $2|x$ by taking (1.1) modulo $p.$

(2)\ Suppose that $p>2$ is a prime and $p\mid gcd(m, n^{2}+1).$
Hence $$a\equiv -n^{r}\equiv -(n^{2})^{r/2}\equiv 1(mod p), b\equiv 0(mod p), c=m^{2}+n^{2}\equiv -1(mod p).$$
It is easy to see that $2|z$ by taking (1.1) modulo $p.$

(3)\ Suppose that $p>2$ is a prime and $p\mid gcd(m^{2}+1, n).$
Hence $$a\equiv m^{r}\equiv (m^{2})^{r/2}\equiv -1(mod p), b\equiv 0(mod p), c=m^{2}+n^{2}\equiv -1(mod p).$$
It is easy to see that $x\equiv z(mod 2)$ by taking (1.1) modulo $p.$

This completes the proof of Lemma 2.1.

\begin{lemma}\label{lemma1}
Suppose that (1.2) holds. If $m\equiv 3(mod 4),$ then $2\mid x.$
\end{lemma}

Proof. By $m\equiv 3(mod 4),$ there is a $d$ with $d\mid m, d\equiv 3(mod 4).$ Then
$$a\equiv -n^{r}\equiv -(n^{2})^{r/2}(mod d), b\equiv 0(mod d), c\equiv n^{2}(mod d).$$
Hence $$(\frac{-(n^{2})^{r/2}}{d})^{x}=(\frac{(n^{2})}{d})^{z},$$ where $(\frac{k}{l})$ is Jacobi symbol.
That is $$(\frac{-1}{d})^{x}=1.$$ Since $(\frac{-1}{d})=(-1)^{(d-1)/2}=-1,$ then $(-1)^{x}=1.$ Hence $2|x.$

This completes the proof of Lemma 2.2.

\begin{lemma}\label{lemma1}~\cite{Mi18}
Suppose that (1.2) holds.

(1)\ Assume that $r\equiv 2(mod 8).$

(i)If $m+n\equiv 7(mod 8),$ then $2\mid y.$

(ii)If $m+n\equiv 3(mod 8),$ then $2\mid z.$

(iii)If $m+n\equiv 5(mod 8)$ or $m-n\equiv \pm3(mod 8),$ then $y\equiv z(mod 2).$

(2)\ Assume that $r\equiv 6(mod 8).$

(i)If $m-n\equiv 7(mod 8),$ then $2\mid y.$

(ii)If $m-n\equiv 3(mod 8),$ then $2\mid z.$

(iii)If $m-n\equiv 5(mod 8)$ or $m+n\equiv \pm3(mod 8),$ then $y\equiv z(mod 2).$
\end{lemma}

\begin{lemma}\label{lemma1}
Suppose that (1.2) holds. If $n\equiv 0(mod 4), m\equiv 3(mod 4),$ then $2\mid y.$
\end{lemma}

Proof. Since $m+n\equiv 3(mod 4),$ then $m-n=m+n-2n\equiv m+n\equiv 3\quad or\quad 7(mod 8).$

(i)$r\equiv 2(mod 8).$ If $m+n\equiv 7(mod 8),$ then $2\mid y$ by Lemma 2.3(1).
If $m+n\equiv 3(mod 8),$ then $2\mid z$ by Lemma 2.3(1). Furthermore $m-n\equiv m+n-2n\equiv m+n\equiv 3(mod 8).$ Hence $y\equiv z(mod 2)$ by Lemma 2.3(1). That is $2\mid y.$

(ii)$r\equiv 6(mod 8).$ If $m-n\equiv 7(mod 8),$ then $2\mid y$ by by Lemma 2.3(2).
If $m-n\equiv 3(mod 8),$ then $2\mid z$ by Lemma 2.3(2). Furthermore $m+n\equiv m-n+2n\equiv m-n\equiv3(mod 8).$ Hence $y\equiv z(mod 2)$ by Lemma 2.3(2). That is $2\mid y.$

This completes the proof of Lemma 2.4.

We may assume that $n>1.$ Then we can define integers $\alpha, \beta(\alpha\geq 1, \beta\geq 2),$ and positive odd integers $i, j$ as follows:
$$m=2^{\beta}j+e, n=2^{\alpha}i, e\in\{1, -1\}, if 2\mid n.$$
$$m=2^{\alpha}i, n=2^{\beta}j+e, e\in\{1, -1\}, if 2\nmid n.$$
Then $\alpha\geq1, \beta\geq 2.$

\begin{lemma}\label{lemma1}~\cite{Mi18}

(1)\ If $2\alpha\neq\beta+1, r\equiv 2(mod 4), y>1,$ then $x\equiv z(mod 2).$

(2)\ If (1.2) holds, $2\alpha=\beta+1, j\equiv e(mod 4),$ then $2\mid z.$
\end{lemma}

If $2\mid x, 2\mid y,$ then we may write that $x=2X, y=2Y.$
\begin{lemma}\label{lemma1}~\cite{Mi18}
If $2\mid x, 2\mid y.$ Then $2\nmid XY.$
\end{lemma}

\begin{lemma}\label{lemma1}
Let $N\in Z^{+}, N>1, A, B, C\in Z,$ gcd$(A, B)=1, AB\neq 0$ we have

(1)~\cite{CD02} If $2\mid A,$ then the equation $A^{2N}+B^{2}=C^{4}$ has no solution $(A, B).$

(2)~\cite{CD02} $A^{2N}+B^{4}=C^{2}$ has no solution $(A, B).$

(3)~\cite{BeCh12, BeSk04, BrCh03} If $N\geq 3,$ then the equation $A^{2}+B^{6}=C^{N}$ has no solution $(A, B).$

(4)~\cite{BeSk04, BrCh03} If $N\geq 4, N\neq 5,$ then the equation $A^{2}+B^{N}=C^{4}$ has no solution $(A, B).$

(5)~\cite{BeChDa15} If $A, B, C$ are nonzero coprime integers for which $A^{2}+B^{2N}=C^{5},$ where $N\geq2$ is an integer, then $B\equiv 1(mod 2).$
\end{lemma}

If $2\mid z,$ then we may write that $z=2Z.$
\begin{lemma}\label{lemma1}
If $r\equiv 2(mod 4), 2\mid x, 2\mid z, 2\mid y.$ Then $2\nmid XYZ.$
\end{lemma}

Proof. By Lemma 2.6, $2\nmid XY.$

In view of $$a^{2}+b^{2}=c^{r}, a^{2X}+b^{2Y}=c^{2Z},$$ $$\max\{a, b\}<c^{r/2}<\min\{a^{2}, b^{2}\}, \max\{a^{X}, b^{Y}\}<c^{Z}<\min\{a^{2X}, b^{2Y}\}.$$
Hence $c^{Z}<a^{2X}<c^{rX}, c^{Z}<b^{2Y}<c^{rY},$ that is $Z<\min\{rX, rY\}.$

If $2\mid Z,$ that is $4\mid z.$ Then $$a^{2X}+b^{2Y}=(c^{z/4})^{4}.$$ By lemma 2.7(1), we have $Y=1.$ In view of lemma 2.7(4), we have $X=1,$ hence $Z=r/2\equiv 1(mod 2),$ a contradiction.

This completes the proof of Lemma 2.8. \hfill$\Box$\\

\section{Estimates of $Y, \Delta$ and $K$}
$ $

We always assume that $(x, y, z)$ is a solution of the equation (1.1). In this section, we will give estimates of $Y, \Delta$ and $K.$
Let $$\Delta=\mid \frac{r}{2}x-z\mid, K=\frac{m}{n}, F_{r}(K)=\frac{c^{r/2}}{a}.$$

For an algebraic number $\gamma$ of degree $d$ over $Q,$ we define as usual the absolute logarithmic height of $\gamma$ by the formula $$h(\gamma)=\frac{1}{d}(\log|a_{0}|+\sum_{i=1}^{d}\log\max\{1, |\gamma^{(i)}|\}),$$
where $a_{0}$ is the leading coefficient of the minimal polynomial of $\gamma$ over $Z$ and the $\gamma^{(i)}$ are the conjugates of $\gamma$ in the field $C$ of complex numbers.

Let $\gamma_{1}, \gamma_{2}$ be algebraic numbers with $\gamma_{1}\gamma_{2}\neq 0.$  let $\log \gamma_{1}$ and $\log \gamma_{2}$ be any determinations of their logarithms. Consider the linear form $$\Lambda=b_{2}\log\gamma_{2}-b_{1}\log\gamma_{1},$$ where $b_{1}, b_{2}\in Z^{+}.$ Without loss of generality, we suppose that $\mid\gamma_{1}\mid, \mid\gamma_{2}\mid\geq 1.$

\begin{lemma}~\cite{Ma00} Corollary 2.3,

Suppose that $\Lambda\neq 0, \gamma_{i}\in C.$
$$A_{i}\geq\max\{Dh(\gamma_{i}), \mid\log\gamma_{i}\mid, 0.16\}, i=1, 2, \quad B=\max\{\mid b_{i}\mid, i=1, 2\}, \quad D=[Q(\gamma_{1}, \gamma_{2}): Q].$$
Then $$\log|\Lambda|>-2^{2.5}e^{2}30^{5}D^{2}A_{1}A_{2}\log(eD)\log(eB).$$
\end{lemma}

\begin{lemma}~\cite{Mi18}
Let $2|x, 2|y.$

(1)\ If $2|z,$ then $$Y\leq\frac{\log4(c-1)}{2(\alpha+1)\log2}.$$
Furthermore if $\exists p\mid n, p\in \{2, 3, 5\},$
then $Y=1$ or $$Y\leq\frac{\log\delta^{2}_{1}(c-1)}{2\log\prod_{q\in S_{1}}(r_{(q)}\cdot n_{(q)})},$$ where $S_{1}=\{2, 3, 5\}\cap\{p: p|n\}$ and
\  \begin{displaymath}\delta_{1}=\begin{cases}
    \ 1, &  \text{$2\nmid n$}, \\
    \ 2, &  \text{$2\mid n$}, \\
 \end{cases}                \end{displaymath}

(2)\ If $2\nmid z,$ then $$Y\leq\frac{\log(c-1)}{2(\alpha+1)\log2}.$$
Furthermore suppose that $\sqrt{c}/\log c>r, Y>1.$ Then $$X\leq\frac{\log(c-4)}{\log {\AA}} \quad Y\leq\frac{\log\delta^{2}_{2}(c-1)}{2\log2^{\alpha+1}\prod_{q\in S_{2}}(r_{(q)}\cdot n_{(q)})}.$$
where $${\AA}=\min\{a_{(p)}\}, S_{2}=\{3, 5\}\cap\{p: p|n\}, \delta_{2}=\prod_{q\in S_{2}}z_{(q)},$$
moreover $z_{(3)}\leq 3,$ and $z_{(5)}=1$ provided $z_{(3)}=3.$
\end{lemma}

If $2|x, 2|y,$ then $$a^{2X}+b^{2Y}=c^{z}.$$ That is $$(a^{X}+b^{Y}\sqrt{-1})(a^{X}-b^{Y}\sqrt{-1})=c^{z}.$$
Since $2\nmid c,$ then two factors on the left-hand side of the above equality are relatively prime in $Z[i].$
Hence we can write $$a^{X}+b^{Y}\sqrt{-1}=(a_{1}+b_{1}\sqrt{-1})^{z},$$ where $a^{2}_{1}+b^{2}_{1}=c.$
Then we can write $$a^{X}+b^{Y}\sqrt{-1}=\epsilon\omega^{z}, \epsilon\in\{\pm 1, \pm i\}, \omega\in Z[i].$$
Without loss of generality, we may write $\omega=c^{1/2}e^{i(\psi/2)}, \psi\in(0, \pi).$
\begin{lemma}
If $2|x, 2|y,$ then $$\min\{a^{X}, b^{Y}\}\geq \pi^{-1}c^{z/2}\mid z\psi-j\pi\mid,$$ where $j\in Z$ such that $$\mid z\psi-j\pi\mid=\min_{j^{\prime}\in Z}\mid z\psi-j^{\prime}\pi\mid.$$
\end{lemma}

Proof. Noting that $$a^{X}+b^{Y}\sqrt{-1}=\epsilon\omega^{z}, \epsilon\in\{\pm 1, \pm i\}, \omega=c^{1/2}e^{i(\psi/2)}\in Z[i].$$
then $$\{a^{X}, b^{Y}\}=\{c^{z/2}\mid\cos(z\psi/2)\mid, c^{z/2}\mid\sin(z\psi/2)\mid\} \eqno{(3.1)}$$

If $0<x<\frac{\pi}{4},$ then $$\sin x>x-\frac{x^{3}}{6}>\frac{2}{\pi}x, \quad \cos x>1-\frac{x^{2}}{2}>\frac{2}{\pi}x.$$
Hence $$\min\{\sin x, \cos x\}>\frac{2}{\pi}x.$$
Noting that $$\mid z\psi/2-j\pi/2\mid<\frac{\pi}{4}.$$
Then $$\min\{\mid\cos(z\psi/2)\mid, \mid\sin(z\psi/2)\mid\}=\min\{\mid\cos(z\psi/2-j\pi/2)\mid, \mid\sin(z\psi/2-j\pi/2)\mid\}>\frac{2}{\pi}\mid z\psi/2-j\pi/2\mid.$$
Hence $$\min\{a^{X}, b^{Y}\}=\min\{c^{z/2}\mid\cos(z\psi/2)\mid, c^{z/2}\mid\sin(z\psi/2)\mid\}\geq\frac{2}{\pi}c^{z/2}\mid z\psi/2-j\pi/2\mid.$$
So $$\min\{a^{X}, b^{Y}\}\geq \pi^{-1}c^{z/2}\mid z\psi-j\pi\mid.$$

This completes the proof of Lemma 3.3. \hfill$\Box$\\

\begin{lemma}
Let $2\mid x, 2\mid y.$ Then $\min\{a^{X}, b^{Y}\}>c^{z/2-6.5\times10^{10}\log z}.$
\end{lemma}

Proof. By Lemma 3.3, $$\min\{a^{X}, b^{Y}\}\geq \pi^{-1}c^{z/2}\mid z\psi-j\pi\mid,$$
where $j\in Z$ such that $$\mid z\psi-j\pi\mid=\min_{j^{\prime}\in Z}\mid z\psi-j^{\prime}\pi\mid.$$
It is clear that $j\geq 0.$ In view of $\mid z\psi-j\pi\mid\leq\frac{\pi}{2},$
$$j\leq z\psi/\pi+\frac{1}{2}<z+\frac{1}{2}.$$ Hence $0\leq j\leq z.$
Let $$(\gamma_{1}, \gamma_{2})=(w, -1), w=\omega/\omega^{\prime}, (b_{1}, b_{2})=(z, -j), \Lambda=z\log w-j\log(-1).$$
Then $$D=2, h(w)=\frac{1}{2}\log c, |\log w|=\psi, h(-1)=0, |\log(-1)|=\pi, \gamma_{i}\in C.$$
Then $\Lambda=(z\psi-j\pi)i\neq 0$ by (3.1). Furthermore $$\min\{a^{X}, b^{Y}\}>\pi^{-1}c^{z/2}\mid\Lambda\mid.$$
Let $$A_{1}\geq\max\{2h(w), \mid \log w\mid \}=\log c, A_{2}\geq\max\{2h(-1), \mid \log(-1)\mid \}=\pi.$$
$$B=\max\{\mid z\mid, \mid -j\mid\}=z.$$
By Lemma 3.1, $$\log|\Lambda|>-2^{2.5}e^{2}30^{5}D^{2}A_{1}A_{2}\log(eD)\log(eB)>-64832990895(\log c)\log z.$$

That is $$\min\{a^{X}, b^{Y}\}>\pi^{-1}c^{z/2-64832990895\log z}>c^{z/2-64832990896\log z}.$$

Hence $$\min\{a^{X}, b^{Y}\}>c^{z/2-6.5\times10^{10}\log z}.$$

This completes the proof of Lemma 3.4. \hfill$\Box$\\

\begin{lemma}
Let $2\mid x, 2\mid y.$ If $rX<z.$ Then $\Delta<26\times10^{10}\log z.$
\end{lemma}

Proof. By Lemma 3.3, $$X>\frac{z/2-6.5\times10^{10}\log z}{\log a}\log c,$$
In view of $a^{2}+b^{2}=c^{r}, \max\{a, b\}<c^{r/2}<\min\{a^{2}, b^{2}\}.$ Then $c^{r/4}<a<c^{r/2}$
\begin{eqnarray*}
\Delta=z-rX&<&z-r\frac{z/2-6.5\times10^{10}\log z}{\log a}\log c\\
&=&\frac{6.5\times10^{10}r\log c}{\log a}\log z-\frac{r/2\log c-\log a}{\log a}z\\
&<&\frac{6.5\times10^{10}r\log c}{\log a}\log z\\
&<&26\times10^{10}\log z\\
\end{eqnarray*}

This completes the proof of Lemma 3.5. \hfill$\Box$\\

\begin{lemma}~\cite{Mi18}
Let $2\mid x, 2\mid y.$

(1)\ If $\Delta=0,$ then $(X, Y, z)=(1, 1, r);$ If $\Delta>0,$ then $\Delta>\frac{\log m}{\log n}.$

(2)\ If $rX>z.$ Then
\  \begin{displaymath}\Delta<\begin{cases}
    \ \frac{\log 4c}{\log c}\frac{\log b^{2}}{\log a}\frac{\log F_{r}(K)}{\log 2^{\alpha+1}}, &  \text{$2\mid z$}, \\
    \ \frac{2\log F_{r}(K)}{\log 3}, &  \text{$2\nmid z, \sqrt{c}/\log c>r$}, \\
 \end{cases}                \end{displaymath}

(3)\ If $rX<z,$ then $z<\frac{2\log b}{\log c}Y+\frac{\log 2}{\log c}, z<rY$ and
\  \begin{displaymath}\Delta<\begin{cases}
    \ \frac{\log b}{\log c}Y+\frac{\log 2}{2\log c}-1, &  \text{$2\mid z$}, \\
    \ \frac{2\log b}{\log c}Y+\frac{\log 2}{\log c}-rX, &  \text{$2\nmid z$}, \\
 \end{cases}                \end{displaymath}
\end{lemma}

\begin{rem}\label{thm1}
Suppose $2\mid x, 2\mid y.$ By Lemma 3.6(1), we only need consider the Case $\Delta>0.$
\end{rem}

\begin{lemma}\label{lemma1}
If $m>2rn, r\equiv 2(mod 4),$ then $F_{r}(K)<\frac{K^{2}}{K^{2}-(32/49)r^{2}}.$
\end{lemma}

Proof. (i)$r>2. F_{r}(K)=\frac{c^{r/2}}{a}=\frac{1}{c^{r/2}/a}=\frac{1}{\cos\phi},$ where $\phi=\arctan\frac{b}{a}.$
Since
\begin{eqnarray*}
\frac{b}{a}&=&\frac{rm^{r-1}n-C_{r}^{3}m^{r-3}n^{3}+\cdots-C_{r}^{r-3}m^{3}n^{r-3}+rmn^{r-1}}
{m^{r}-C_{r}^{2}m^{r-2}n^{2}+\cdots+C_{r}^{r-2}m^{2}n^{r-2}+n^{r}}\\
&<&\frac{rm^{r-1}n}{m^{r}-C_{r}^{2}m^{r-2}n^{2}}\\
&=&\frac{r}{K}\frac{K^{2}}{K^{2}-\frac{r(r-1)}{2}}\\
&=&\frac{r}{K}\frac{K^{2}}{K^{2}-\frac{K^{2}}{8}}\\
&<&\frac{8r}{7K}
\end{eqnarray*}
Then $$\phi=\arctan\frac{b}{a}<\arctan\frac{8r}{7K}<\frac{8r}{7K}.$$
It follow that $$\cos\phi>\cos\frac{8r}{7K}>1-\frac{64r^{2}/49K^{2}}{2}=1-\frac{32r^{2}}{49K^{2}}.$$
Hence $$F_{r}(K)=\frac{1}{\cos\phi}<\frac{1}{1-\frac{32r^{2}}{49K^{2}}}=\frac{K^{2}}{K^{2}-(32/49)r^{2}}.$$

(ii)$r=2. F_{r}(K)=F_{2}(K)=\frac{K^{2}+1}{K^{2}-1}=1+\frac{2}{K^{2}-1}<1+\frac{128/49}{K^{2}-128/49}=\frac{K^{2}}{K^{2}-(32/49)r^{2}}.$

This completes the proof of Lemma 3.7. \hfill$\Box$\\

\section{The case $y=1$}
$ $

We always assume that $(x, y, z)$ is a solution of the equation (1.1). In this section, we consider the case $y=1.$ In view of $a^{2}+b^{2}=c^{r}$ and gcd$(a^{2}+b^{2}, a+b)=1, x\geq2.$

\begin{lemma}~\cite{La08} Corollary 2, $m=10$
Suppose that $\gamma_{1}, \gamma_{2}\in R^{+}$ are multiplicatively independent with $\log\gamma_{1}, \log\gamma_{2}\in R^{+}.$
$$A_{i}\geq\max\{h(\gamma_{i}), \mid \log\gamma_{i}\mid/D, 1/D \}, i=1, 2. \quad D=[Q(\gamma_{1}, \gamma_{2}): Q] \quad b^{\prime}=\frac{b_{1}}{DA_{2}}+\frac{b_{2}}{DA_{1}}.$$
Then $$\log|\Lambda|>-25.2D^{4}(\max\{\log b^{\prime}+0.38, \frac{10}{D}, 1\})^{2}A_{1}A_{2}.$$
\end{lemma}

\begin{lemma}\label{lemma1}
Let $y=1, s=\frac{x}{\log c}, \Lambda=z\log c-x\log a,$ then $$\log\Lambda>-25.2(\max\{\log(2s+1)+0.38, 10\})^{2}(\log a)\log c.$$
\end{lemma}

Proof. Since $a^{2}+b^{2}=c^{r},$ then $$\max\{a, b\}<c^{r/2}<\min\{a^{2}, b^{2}\}.$$
In view of $b<a^{2}\leq a^{x}, c^{z}=a^{x}+b<2a^{x}.$ Hence $$\frac{z}{\log a}<\frac{\log 2}{(\log a)\log c}+\frac{x}{\log c}<\frac{x}{\log c}+1.$$
Let $(\gamma_{1}, \gamma_{2})=(a, c), (b_{1}, b_{2})=(x, z).$
Then $$D=1, \log A_{i}=h(\gamma_{i})=\log\gamma_{i}, b^{\prime}=\frac{x}{\log c}+\frac{z}{\log a}<2\frac{x}{\log c}+1=2s+1, s=\frac{x}{\log c}.$$
By Lemma 4.1,
\begin{eqnarray*}
\log\Lambda &>&-25.2D^{4}(\max\{\log b^{\prime}+0.38, \frac{10}{D}, 1\})^{2}(\log A_{1})\log A_{2}\\
&>&-25.2(\max\{\log(2s+1)+0.38, 10\})^{2}(\log a)\log c\\
\end{eqnarray*}

This completes the proof of Lemma 4.2.

\begin{lemma}\label{lemma1}
If $y=1,$ then $$x<7531.1\log c, \Delta<7531.1\log F_{r}(K).$$
\end{lemma}

Proof. Let $\Lambda=z\log c-x\log a,$ then $\Lambda>0.$ Noting that $$a^{2}+b^{2}=c^{r},$$ then $$\max\{a, b\}<c^{r/2}<\min\{a^{2}, b^{2}\}.$$ Hence $$\Lambda<e^{\Lambda}-1=c^{z}/a^{x}-1=b/a^{x}<a^{-(x-2)},$$ that is $\log\Lambda<-(x-2)\log a.$
By Lemma 4.2, $$(x-2)\log a<25.2(\max\{\log(2s+1)+0.38, 10\})^{2}(\log a)\log c.$$
that is $$s=\frac{x}{\log c}<\frac{2}{(\log a)\log c}+25.2(\max\{\log(2s+1)+0.38, 10\})^{2}.$$

If $\log(2s+1)+0.38<10,$ then $s<7531.1$

If $\log(2s+1)+0.38>10,$ then $s<25.2(\log(2s+1)+0.38)^{2}+1,$ then $s<7500,$ a contradiction.

Hence $x/\log c=s<7531.1$ that is $x<7531.1\log c.$

Noting that $a^{2}+b^{2}=c^{r},$ then $\max\{a, b\}<c^{r/2}<\min\{a^{2}, b^{2}\}.$ Hence $c^{r/4}<a<c^{r/2}.$ By $x\geq 2, c^{z}=a^{x}+b<2a^{x}<2c^{rx/2},$ that is $c^{2z}<4c^{rx}<c^{rx+1}.$ Hence $2z\leq rx.$
In view of $$a^{x}<c^{z},$$
We have $$\Delta=\frac{r}{2}x-z<\frac{r}{2}x-\frac{\log a}{\log c}x=\frac{x}{\log c}(\frac{r}{2}\log c-\log a)=\frac{x}{\log c}\log F_{r}(K).$$
Then $$\frac{\Delta}{\log F_{r}(K)}<\frac{x}{\log c} \eqno{(4.1)}$$
Hence $\Delta<7531.1\log F_{r}(K).$

This completes the proof of Lemma 4.3.

\begin{lemma}\label{lemma1}~\cite{Mi18}
If $r=2, y=1,$ then $$x<1534\log c, \Delta<1534\log F_{2}(K)=1534\log\frac{c}{a}, K<56.$$
\end{lemma}

\begin{lemma}\label{lemma1}
Let $p$ be a prime, $U, V$ be distinct coprime integers such that
\  \begin{displaymath}\begin{cases}
    \ U\equiv V(mod p), &  \text{$p>2$}, \\
    \ U\equiv V(mod 4), &  \text{$p=2$}, \\
 \end{cases}                \end{displaymath}
Then $$ord_{p}(U^{k}-V^{k})=ord_{p}(U-V)+ord_{p}(k), k\in Z^{+}.$$
\end{lemma}

Noting that if $r=2,$ then $a^{2}+b^{2}=c^{2},$ hence $$a=m^{2}-n^{2}, b=2mn, c=m^{2}+n^{2},$$ furthermore $$\max\{a, b\}<c<\min\{a^{2}, b^{2}\}.$$
\begin{lemma}\label{lemma1}
If $y=1, r=2, n\geq 4, 2\mid n, m\equiv 3(mod 4),$ then
$$\Delta>\frac{1}{m+1}(\frac{n_{(q)}^{2}}{1534^{2}\log^{2} F_{2}(K)}-n)$$ for some $q\in\{3, 5\}.$
\end{lemma}

Proof. If gcd$(x, z)=d>1,$ then $$2mn=b=c^{z}-a^{x}=(c^{d}-a^{d})[(c^{z/d})^{d-1}+\cdots+(a^{x/d})^{d-1}],$$ which contradicts $2mn<m^{2}+n^{2}=c.$ Hence gcd$(x, z)=d=1.$ By Lemma 2.2, $2\mid x.$ Then $2\nmid z.$ Hence $2\nmid\Delta.$

If $z=1,$ then $$(m+n)^{2}(m-n)^{2}<(m^{2}-n^{2})^{x}=a^{x}=c-b=(m^{2}+n^{2})-2mn=(m-n)^{2},$$ a contradiction. Hence $z\geq 3.$

If $z>x,$ then $$2mn=b=c^{z}-a^{x}>c^{x}-a^{x}=(c-a)(c^{x-1}+\cdots+a^{x-1}),$$ which contradicts $2mn<m^{2}+n^{2}=c.$ Hence $x>z\geq3.$ That is $x\geq 4.$

We consider the case $q|n, q\in\{3, 5\}.$ Then $q\nmid m.$ Taking (1.1) modulo $n^{2},$ we have $$m^{2z}(m^{2\Delta}-1)\equiv-2mn(mod n^{2}).$$
Then $$ord_{q}(m^{2\Delta}-1)=ord_{q}(n) \eqno{(4.2)}.$$
Noting that $m^{2}\equiv \pm 1(mod q), q\in\{3, 5\}.$
In view of $2\nmid \Delta,$ $$m^{2}\equiv m^{2\Delta}\equiv 1(mod q).$$
By Lemma 4.5, $$ord_{q}(m^{2\Delta}-1)=ord_{q}(m^{2}-1)+ord_{q}(\Delta).$$
By (4.2) $$ord_{q}(\Delta)=ord_{q}(m^{2\Delta}-1)-ord_{q}(m^{2}-1)=ord_{q}(n)-ord_{q}(m-\epsilon), \epsilon\in\{1, -1\}.$$
Let $M=m-\epsilon=M_{(q)}T$ with $q\nmid T.$
Then $$\frac{n_{(q)}}{M_{(q)}}\mid \Delta,$$ that is $$\Delta\geq\frac{n_{(q)}}{M_{(q)}}.$$
By Lemma 4.4, $$M_{(q)}>\frac{n_{(q)}}{1534\log F_{2}(K)}.$$
Since $$m^{2}\equiv 1+2\epsilon M_{(q)}T(mod M^{2}_{(q)}), \quad n\equiv 0(mod M_{(q)}),$$
then $$a^{x}\equiv m^{2x}\equiv(1+2\epsilon M_{(q)}T)^{x}\equiv1+2\epsilon M_{(q)}Tx(mod M^{2}_{(q)}),$$
$$b=2mn=2(\epsilon+M_{(q)}T)n\equiv 2\epsilon n(mod M^{2}_{(q)}),$$
$$c^{z}\equiv m^{2z}\equiv1+2\epsilon M_{(q)}Tz(mod M^{2}_{(q)}),$$
By $$a^{x}+b=c^{z},$$ $$M_{(q)}T\Delta+n\equiv 0(mod M^{2}_{(q)}).$$
So $$\Delta\geq\frac{M^{2}_{(q)}-n}{M_{(q)}T}=\frac{M^{2}_{(q)}-n}{m-\epsilon}\geq\frac{M^{2}_{(q)}-n}{m+1}
>\frac{1}{m+1}(\frac{n_{(q)}^{2}}{1534^{2}\log^{2} F_{2}(K)}-n).$$

This completes the proof of Lemma 4.6.

\begin{lemma}\label{lemma1}
If $y=1, r=2, n\geq 4, 2\mid n, m\equiv 3(mod 4)$ and $$n_{(q)}\geq (1534\log n)^{1.5}t(n)\sqrt{n}$$ for some $q\in\{3, 5\}.$
Then the equation has no solution.
Where $$t(n)=\sqrt{\frac{1}{1534}+(1+\frac{4}{n})\log(1+\frac{2n^{2}}{6n+9})}\frac{\log(1+\frac{2n^{2}}{6n+9})}{(\log n)^{1.5}}.$$
\end{lemma}

Proof. By Lemma 4.4, 4.6, $$\frac{1}{m+1}(\frac{n_{(q)}^{2}}{1534^{2}\log^{2}(F_{2}(K))}-n)<1534\log F_{2}(K).$$
that is $$n_{(q)}n^{-1/2}<1534\sqrt{1+1534(K+1/n)\log F_{2}(K)}\log F_{2}(K).$$
Clearly, the right-hand side is a decreasing function of $K.$ By~\cite{De}, we may assume that $m\geq n+3.$ Then $K\geq 1+3/n.$
Hence $$n_{(q)}n^{-1/2}<1534\sqrt{1+1534(1+4/n)\log(1+\frac{2n^{2}}{6n+9})}\log(1+\frac{2n^{2}}{6n+9}).$$
Hence $n_{(q)}<(1534\log n)^{1.5}t(n)\sqrt{n},$ a contradiction.

This completes the proof of Lemma 4.7.

\section{Proof of Theorem 1.1 and Corollary 1.1}
$ $

We always assume that $(x, y, z)$ is a solution of the equation (1.1).
Recall that $$m=2^{\beta}j+e, n=2^{\alpha}i, 2\nmid ij, e\in\{1, -1\}, if 2\mid n.$$
$$m=2^{\alpha}i, n=2^{\beta}j+e, 2\nmid ij, e\in\{1, -1\}, if 2\nmid n.$$ Then $\alpha\geq1, \beta\geq 2.$

\begin{lemma}\label{lemma1}
Suppose that (1.2) holds. If $$n\equiv 2(mod 4), m\equiv 3 (mod 4), r\equiv 2(mod 4), F_{r}(K)>e^{1/7531.1},$$
Then $(x,y,z)=(2,2,r).$
\end{lemma}

Proof. By Lemma 4.3, we have $y>1.$ In view of Lemma 2.2, $2\mid x.$ Since $\alpha=1, \beta\geq 2,$
then $2\alpha\neq\beta+1.$ By lemma 2.5, we have $2\mid z.$
Noting that $m+n\equiv 1(mod 4),$ then $m+n\equiv 1, 5(mod 8).$
If $m+n\equiv 5(mod 8),$ then $y\equiv z(nod 2)$ by lemma 2.3. Hence $2\mid y.$
If $m+n\equiv 1(mod 8),$ then $m-n=m+n-2n\equiv 5(mod 8).$ By lemma 2.3, $y\equiv z(nod 2).$ Hence $2\mid y.$
In any way, we have $2\mid y.$

Then $(a^{X}, b^{Y}, c^{Z})$ forms a primitive Pythagorean triple, hence
$$a^{X}=k^{2}-l^{2}, b^{Y}=2kl, c^{Z}=k^{2}-l^{2}$$ with $k>l\geq 1,$ gcd$(k, l)=1, k-l\equiv 1(mod 2).$
Since$$a^{2}+b^{2}=c^{r}, a^{2X}+b^{2Y}=c^{2Z},$$
then $$\max\{a, b\}\leq c^{r/2}<\min\{a^{2}, b^{2}\}\quad \quad\max\{a^{X}, b^{Y}\}\leq c^{Z}<\min\{a^{2X}, b^{2Y}\}.$$
Hence $a^{X}<b^{2Y}<a^{4Y}.$ Then $X<4Y.$

We have $$a\equiv m^{r}-\frac{r(r-1)}{2}m^{r-2}n^{2}\equiv1-4\equiv 5(mod 8).$$
It follow from $2\nmid X$ that $k^{2}-l^{2}=a^{X}\equiv 5(mod 8),$ that is $l\equiv 2(mod 4), 2\nmid k.$
Hence $$Y=ord_{2}(2kl)/ord_{2}(b)=2/ord_{2}(2n)=1.$$
We have $X<4.$ Suppose that $X=3,$ then $Z=1$ by lemma 2.7(3). That is $a^{6}+b^{2}=c^{2},$ which contradict $a^{2}+b^{2}=c^{r}.$ Therefore $X=1,$ then $Z=r/2.$ Hence $(x,y,z)=(2,2,r).$

This completes the proof of Lemma 5.1. \hfill$\Box$\\

\begin{rem}\label{thm1}
In fact, during the proof of Lemma 2, condition $F_{r}(K)>e^{1/7531.1}$ is only needed when discussing the case $y=1.$
\end{rem}

\begin{lemma}\label{lemma1}
If $x>2t\log t, t>10^{10},$ then $x>t\log x.$
\end{lemma}

Proof. Let $f(x)=x-t\log x,$ then $f^{\prime}(x)=1-t/x>0.$ Hence $f(x)$ is a increasing function.
Noting that $$f(2t\log t)=2t\log t-t\log(2t\log t)=t(\log t-\log 2-\log\log t).$$
Let $$g(t)=\log t-\log 2-\log\log t,$$
then $$g^{\prime}(t)=\frac{1}{t}-\frac{1}{t\log t}>0.$$
Hence $g(x)$ is a increasing function.
Noting $g(10^{10})>0.$ Then $g(t)>0.$
Hence $f(2t\log t)>0.$ Then $f(x)>0.$
So $x>t\log x.$

This completes the proof of Lemma 5.2. \hfill$\Box$\\

\textbf{Proof of Theorem 1.1.} If $y=1,$ by Lemma 3.7, 4.3, $$\frac{K^{2}}{K^{2}-\frac{32}{49}r^{2}}>F_{r}(K)>e^{1/7531.1}.$$ Then  $K<70.2r,$ that is $m<70.2nr,$ a contradiction. Hence $y>1.$ Then we can assume that $n\equiv 0(mod 4)$ by lemma 5.1. That is $\alpha\geq 2.$ Hence $2\mid x, 2\mid y$ by Lemma 2.2, 2.4. Suppose $\Delta>0.$ In view of $m>3e^{r},$ $$\frac{\sqrt{c}}{\log c}>\frac{m}{\log m^{3}}=\frac{3e^{r}}{3\log 3e^{r}}>r.$$

(i)$rX>z.$ By lemma 3.6(2),
\  \begin{displaymath}1\leq\Delta<\begin{cases}
    \ \frac{4\log 4c}{(\log 8)\log c}\log F_{r}(K),  &  \text{$2\mid z$}, \\
    \ \frac{2}{\log 3}\log F_{r}(K),  &  \text{$2\nmid z$}, \\
 \end{cases}                \end{displaymath}
Then $$\frac{1}{\log F_{r}(K)}<\max\{\frac{4\log 4c}{(\log 8)\log c}, \frac{2}{\log 3}\}=\max\{\frac{8}{3\log c}+\frac{4}{3\log 2}, \frac{2}{\log 3}\}=2.57.$$
By Lemma 4.1, $\frac{K^{2}}{K^{2}-r^{2}}>F_{r}(K)>e^{\frac{1}{2.57}},$ then $K<3r,$ that is $m<3nr$ contadict $m>70.2nr.$

(ii)$rX<z,$ then $z<rY$ by lemma 3.6(3). Noting that $n\geq4, m\geq 7, \alpha\geq 2.$ It follows from lemma 3.2 that
$$Y\leq\frac{\log(4c)}{2(\alpha+1)\log2}=\frac{\log(4n^{2}+4m^{2})}{2(\alpha+1)\log2}<\frac{\log((m-4)m^{2}+4m^{2})}{2(\alpha+1)\log2}=\frac{\log(m^{3})}{2(\alpha+1)\log2}<\log m.$$
In view of Lemma 3.5, 3.6(1), we have $$\frac{\log m}{\log n}<\Delta<26\times10^{10}\log z<26\times10^{10}\log (r\log m).$$
By $m>e^{r},$ $$\log m<26\times10^{10}(\log n)\log (r\log m)<5.2\times10^{11}(\log n)\log(\log m).$$

On the other hand, $$m>n^{10.4\times10^{11}\log(5.2\times10^{11}\log n)}=e^{10.4\times10^{11}(\log n)\log(5.2\times10^{11}\log n)},$$
then $$\log m>10.4\times10^{11}(\log n)\log(5.2\times10^{11}\log n).$$
Let $t=5.2\times10^{11}\log n,$ then $\log m>2t\log t.$ By Lemma 5.2, $\log m>t\log\log m.$
that is $$\log m>5.2\times10^{11}(\log n)\log(\log m),$$ a contradiction.

This completes the proof of Theorem 1.1.

Proof of Corallary 1.1. Since $r=2, m>n^{10^{11}}>3^{6}n>70.2nr>3^{3}>3e^{2},$
then $$m>\max\{n^{10.4\times10^{11}\log(5.2\times10^{11}\log n)}, 3e^{r}, 70.2nr\}=n^{10.4\times10^{11}\log(5.2\times10^{11}\log n)}.$$

This completes the proof of Corallary 1.1.

\section{Proof of Theorem 1.2 and Theorem 1.3}
$ $

We always assume that $(x, y, z)$ is a solution of the equation (1.1). In this section, we always assume that $r=2, n\geq 4.$ Then $a^{2}+b^{2}=c^{2},$ Hence $$a=m^{2}-n^{2}, b=2mn, c=m^{2}+n^{2},$$ furthermore $$\max\{a, b\}<c<\min\{a^{2}, b^{2}\}.$$
Recall that $$m=2^{\beta}j+e, n=2^{\alpha}i, 2\nmid ij, e\in\{1, -1\}, if 2\mid n.$$
$$m=2^{\alpha}i, n=2^{\beta}j+e, 2\nmid ij, e\in\{1, -1\}, if 2\nmid n.$$ Then $\alpha\geq1, \beta\geq 2.$
If $y>1.$ In view of the proof of Lemma 5.1, we may assume that $\alpha\geq 2.$ By [Mi09, Theorem 1.5]: If $m\equiv4(mod 8) and n\equiv7(mod 16),$ or $m\equiv7(mod 16) and n\equiv4(mod 8),$ then Conjecture 1.1 holds.
We only need to consider the case $2\alpha\neq \beta+1$ when $\alpha=2.$

Let $m\equiv 3(mod 4).$ If $(m, n)$ are not the follow case $2\alpha=\beta+1, j\equiv 1(mod 4), \alpha\geq 3,$ then $x\equiv z(mod 2)$ by Lemma 2.5.
For example $$2\mid\beta or 2\alpha\neq\beta+1 or 2\alpha=\beta+1, j\equiv 3(mod 4).$$

\begin{lemma}\label{lemma1}~\cite{Mi18}
Suppose $K=m/n>1.5, 2\mid x, 2\mid y, \Delta>0.$ Then $x<z<y.$
\end{lemma}

\begin{lemma}\label{lemma1}~\cite{Mi18}
Suppose $K=m/n>1.5, 2\mid x, 2\mid y, \Delta>0.$ Then

(i)gcd$(XY, 6)=1, t\nmid z, t\in \{9, 15\}.$

(ii)Suppose that $2\mid z,$ then $\min\{X, Y\}\geq 5,$ gcd$(Z, 30)=1.$ Moreover, if $K>1.5,$ then $Y\geq 11.$

(iii)Suppose that $2\nmid z, K>1.5,$ then $Y\geq 5, 2Y-3\geq z>Y.$
\end{lemma}

\textbf{Proof of Theorem 1.2.} By Lemma 4.4, we may assume that $y>1.$ In view of Lemma 2.2 and 2.4, we have $2|x, 2|y.$

(1)\ $(m, n)$ are not the case $2\alpha=\beta+1, j\equiv 1(mod 4), \alpha\geq 3.$ Then $2\mid z.$
By Lemma 3.2, $$Y<\frac{\log (4c)}{\log 4^{\alpha+1}} \eqno{(6.1)}$$
Suppose $\Delta>0.$ Since $m>56n,$ then $x<z$ by Lemma 6.1.
In view of 3.6, $$\frac{\log m}{\log n}<\Delta<\frac{\log b}{\log c}\frac{\log(4c)}{\log 4^{\alpha+1}}+\frac{\log 2}{2\log c}-1.$$
Then $$\frac{\log(mn)}{\log n}<\frac{\log b}{\log 4^{\alpha+1}}\frac{\log(4c)}{\log c}+\frac{\log 2}{2\log c}.$$
That is $$\frac{\log4^{\alpha+1}}{\log b}\frac{\log(mn)}{\log n}<\frac{\log(4c)}{\log c}+\frac{\log4^{\alpha+1}}{\log b}\frac{\log 2}{2\log c}.$$
Let $$\nu=\frac{\log b}{\log(mn)}=1+\frac{\log 2}{\log(mn)}.$$
Then $$\frac{\log4^{\alpha+1}}{\nu\log n}<\frac{\log(4c)}{\log c}+\frac{(\log 2)\log4^{\alpha+1}}{2\nu(\log c)\log(mn)}.$$
By Lemma 6.2(ii), $Y\geq 11.$ Then $4c>4^{11(\alpha+1)}$ by (6.1).
Hence $$mn\geq4m>\sqrt{4c}>2^{11(\alpha+1)}.$$
Then $$\frac{\log4^{\alpha+1}}{\nu\log n}<\frac{\log(4c)}{\log c}+\frac{\log4}{22\nu(\log c)}.$$
That is $$\frac{\log(4^{\alpha+1}/n^{\nu})}{\nu\log n}<\frac{\log 4}{\log c}+\frac{\log4}{22\nu\log c}.$$
Hence $$\frac{\log(4^{\alpha+1}/n^{\nu})}{\nu(\log 4)\log n}<\frac{22\nu+1}{22\nu\log c}$$
Noting that $n^{\nu}<4^{\alpha+1}.$
So $$\log c<\frac{(22\nu+1)(\log 2)}{11\log(4^{\alpha+1}/n^{\nu})}\log n=\eta\log n, \quad \eta=\frac{(22\nu+1)(\log 2)}{11\log(4^{\alpha+1}/n^{\nu})} \eqno{(6.2)}$$
which contradicts $m>\sqrt{n^{\eta}-n^{2}}.$

(2)\ $n=2^{\alpha}.$ If $2\mid z.$ Since $n^{\nu}<4^{\alpha+1},$ it is seem to (1), we have (6.2).
In view of $\nu=1+\frac{\log2}{\log(mn)}<1+\frac{\log2}{\log(2^{2\alpha})}\leq 1.25, \eta<2,$ contradicts (6.2).
If $2\nmid z.$ by Lemma 3.2, $$Y<\frac{\log c}{\log 4^{\alpha+1}}.$$
Suppose $\Delta>0.$ In view of $m>56n>1.5n,$ $x<z$ by Lemma 6.1. It follow from Lemma 3.6 that
$$\frac{\log m}{\log n}<\Delta<\frac{2\log b}{\log c}\frac{\log c}{\log 4^{\alpha+1}}+\frac{\log 2}{\log c}-2.$$
Then $$\frac{\log(mn^{2})}{\log n}<\frac{\log b}{\log 2^{\alpha+1}}+\frac{\log 2}{\log c}.$$
Hence $$\frac{(\log (mn^{2}))\log 2^{\alpha+1}-(\log b)\log n}{(\log n)\log 2^{\alpha+1}}<\frac{\log 2}{\log c}.$$

In view of $n=2^{\alpha},$ we have $$\frac{(\log 2)\log m+\alpha(\alpha+1)\log^{2}2}{(\log n)\log 2^{\alpha+1}}<\frac{\log 2}{\log c}.$$

Hence $$\log c<\frac{(\alpha+1)\log 2}{\log m+\alpha(\alpha+1)\log 2}\log n,$$
Noting that $$\frac{(\alpha+1)\log 2}{(\alpha+1)\alpha\log 2+\log m}<1,$$ a contradiction.

This completes the proof of Theorem 1.2.

\textbf{Proof of Theorem 1.3.} By lemma 4.4, 4.7, we have $y>1.$ By Lemma 5.1, we may assume that $\alpha\geq 2.$ In view of Lemma 2.2, 2.4, $2|x, 2|y.$

(1)\ Since $m\equiv 3(mod 4)$ and $(m, n)$ not the case $$2\alpha=\beta+1, j\equiv 1(mod 4), \alpha\geq 3.$$
Then $2\mid z$ by Lemma 2.5. It follow from Lemma 3.2 that
$$Y<\frac{\log 4(c-1)}{2\log n^{\frac{1}{2}+\kappa}} \eqno{(6.2)}$$ with $n^{\frac{1}{2}+\kappa}=2n_{(2)}n_{(3)}n_{(5)}.$ In view of the assumption, we have $\kappa>0.$
Since $Y\geq5$ by Lemma 6.2, it follow that $$n^{5}<n^{(1+2\kappa)Y}<4c=4m^{2}+4n^{2}.$$ Then $m>3n.$ Hence $x<z$ by Lemma 6.1.
In view of Lemma 3.6, we have
$$\frac{\log m}{\log n}<\Delta<\frac{\log b}{\log c}\frac{\log 4c}{(1+2\kappa)\log n}+\frac{\log2}{2\log c}-1.$$
Then $$\frac{\log(mn)}{\log n}<\frac{\log b}{\log c}\frac{\log 4c}{(1+2\kappa)\log n}+\frac{\log2}{2\log c}.$$
That is $$\frac{2(\log c)\log(mn)-(\log2)\log n}{2(\log c)\log n}<\frac{\log b}{\log c}\frac{\log 4c}{(1+2\kappa)\log n}.$$
Hence $$\frac{2(\log c)\log(mn)-(\log2)\log n}{2}<\frac{(\log b)\log 4c}{(1+2\kappa)}.$$
That is $$1+2\kappa<\frac{2(\log b)\log 4c}{2(\log c)\log(mn)-(\log2)\log n}.$$
Then
\begin{eqnarray*}
\kappa &<&\frac{2(\log b)\log 4c-2(\log c)\log(mn)+(\log2)\log n}{4(\log c)\log(mn)-2(\log2)\log n}\\
&=&\frac{2(\log2+\log(mn))(\log4+\log c)-2(\log c)\log(mn)+(\log2)\log n}{4(\log c)\log(mn)-2(\log2)\log n}\\
&=&\frac{2(\log4)\log(mn)+2(\log2)\log c+2(\log2)\log4+(\log2)\log n}{4(\log c)\log(mn)-2(\log2)\log n}\\
&=&\frac{4(\log2)\log m+5(\log2)\log n+2(\log2)\log c+4\log^{2}2}{4(\log c)\log(mn)-2(\log2)\log n}\\
&=&\frac{4\log m+5\log n+2\log c+4\log2}{4(\log c)(\log(mn)/\log n)-2\log2}\cdot\frac{\log2}{\log n}\\
\end{eqnarray*}
Noting that $$f(m)=\frac{4\log m+5\log n+2\log c+4\log2}{4(\log c)(\log(mn)/\log n)-2\log2}=\frac{2\log4m^{2}+5\log n+2\log c}{4(\log c)\log(mn)-2(\log2)\log n}\log n$$ is a decreasing function of $m.$
In view of Lemma 6.2, $Y\geq 11.$ By (6.2), $$4m^{2}+4n^{2}=4c>n^{(1+2\kappa)Y}>n^{11}.$$
Then $$f(m)\leq\frac{2\log(n^{11}-4n^{2})+5\log n+2\log(n^{11}/4)}{2(\log(n^{11}/4))\log((n^{11}-4n^{2})n^{2}/4)-2(\log2)\log n}\log n=s(n)$$
that is $$\kappa<\frac{s(n)\log2}{\log n},$$
Hence $2n_{(2)}n_{(3)}n_{(5)}=n^{1/2+\kappa}<2^{s(n)}\sqrt{n},$ which contradicts the assumption.

(2)\ $(m, n)$ are the case $2\alpha=\beta+1, j\equiv 1(mod 4).$ If $2\mid z.$ It is seem to (1). If $2\nmid z.$ Then $\Delta>0, 2\nmid\Delta.$
We divide the argument into two cases:

Case 1: $x>z.$ Then $m\leq1.5n$ by Lemma 6.1. In view of the assumption, we have $$n_{(q)}>2\sqrt{n\log n}$$ for some $q\in\{3, 5\}.$
Taking the equation (1.1) modulo $n^{2}, m^{2\Delta}\equiv 1(mod n^{2}),$ then $$ord_{q}(m^{2\Delta}-1)\geq ord_{q}(n^{2}).$$
Noting that $$m^{2}\equiv \pm 1(mod q), q\in\{3, 5\}, 2\nmid \Delta.$$
Then $$m^{2}\equiv m^{2\Delta}\equiv 1(mod q),$$
By Lemma 4.5, $$ord_{q}(m^{2}-1)+ord_{q}(\Delta)=ord_{q}(m^{2\Delta}-1)>ord_{q}(n^{2}).$$
Then $$(m+\epsilon)\Delta\equiv0(mod n_{(q)}^{2}), \epsilon\in\{\pm 1\}.$$
Then $$4n\log n<n_{(q)}^{2}\leq (m+1)\Delta \eqno{(6.3)}$$
By $m\geq n+1,$ $$F_{2}(K)=\frac{c}{a}=1+\frac{2n^{2}}{m^{2}-n^{2}}<1+\frac{2n^{2}}{2n+1}=\frac{2n^{2}+2n+1}{2n+1}.$$
From Lemma 3.6(2), $$\Delta<\frac{2}{\log3}\log F_{2}(K)<\frac{2}{\log3}\log\frac{2n^{2}+2n+1}{2n+1}<2\log n.$$
Hence $$4n\log n<n_{(q)}^{2}\leq(m+1)\Delta<(3n+2)\log n,$$ a contradiction.

Case 2: $x<z.$ In view of Lemma 3.2, $$Y<\frac{\log c}{2\log 2^{\alpha+1}}.$$ If $m\leq 1.5n.$ By lemma 3.6, $z<y$ and
\begin{eqnarray*}
\Delta &<&\frac{2\log b}{\log c}Y+\frac{\log 2}{\log c}-rX\\
&<&\frac{2\log b}{\log c}\frac{\log c}{2\log 2^{\alpha+1}}+\frac{\log 2}{\log c}-rX\\
&=&\frac{\log(2mn)}{\log 2^{\alpha+1}}+\frac{\log 2}{\log c}-rX\\
&<&\frac{\log(3n^{2})}{\log 2^{\alpha+1}}-1\\
&=&\frac{\log(3n^{2}/2^{\alpha+1})}{\log 2^{\alpha+1}}<\frac{\log n^{2}}{\log 2^{3}}<2\log n.\\
\end{eqnarray*}
It is seem to case 1, we get a contradicion. Hence $m>1.5n.$ By Lemma 6.1, $x<z<y.$ In view of $2\mid x, 2\mid y, y\geq 4.$

If $5\mid z,$ in view of Lemma 2.7(5), $b=2mn\equiv 1(mod 2),$ a contradiction. Hence $5\nmid z.$
By Lemma 3.2, $$Y<\frac{\log c+2\log\delta_{2}}{2\log Q}.$$
with $$\delta_{2}=\prod_{q\in S_{2}}z_{(q)}=z_{(3)}, Q=2^{\alpha+1}\prod_{q\in S_{2}}n_{(q)}, S_{2}=\{3, 5\}\cap\{p: p\mid n\}.$$
By Lemma 6.2(1), $\delta_{2}=z_{(3)}\leq 3<Q.$
Then $$2Y\log Q-2Y\log\delta^{1/Y}_{2}<\log c.$$ That is $$Y<\frac{\log c}{2\log(Q/\delta_{2}^{1/Y})}.$$
By Lemma 3.6,
$$\frac{\log m}{\log n}<\Delta<\frac{2\log b}{\log c}Y+\frac{\log 2}{\log c}-rX<\frac{\log b}{\log(Q/\delta_{2}^{1/Y})}+\frac{\log 2}{\log c}-2.$$
Then $$\frac{\log(mn^{2})}{\log n}<\frac{\log b}{\log(Q/\delta_{2}^{1/Y})}+\frac{\log 2}{\log c},$$
That is $$\frac{(\log c)\log(mn^{2})-(\log 2)\log n}{(\log c)\log n}<\frac{\log b}{\log(Q/\delta_{2}^{1/Y})},$$
By $n\geq 4,$
\begin{eqnarray*}
\log(Q/\delta_{2}^{1/Y})&<&\frac{(\log b)\log c}{(\log c)\log(mn^{2})-(\log2)\log n}\log n\\
&<&\frac{(\log b)\log c}{(\log c)\log(mn^{2})-\log n}\log n\\
&=&\frac{\log b}{\log(mn^{2})-(\log n)/(\log c)}\log n\\
&<&\frac{\log b}{\log(mn^{2})-1/2}\log n\\
&=&\frac{\log (2mn)}{\log[(n/\sqrt{e})mn]}\log n<\log n\\
\end{eqnarray*}
Then $\delta_{2}^{1/Y}>Q/n.$ By the assumption, $Q=2^{\alpha+1}\prod_{q\in \{3, 5\}}n_{(q)}>3^{1/12}n,$ then $\delta_{2}^{1/Y}>Q/n>3^{1/12}>1.$ Then $\delta_{2}>1.$ Hence $\delta_{2}=z_{(3)}=3.$
By Lemma 6.2(iii), $2Y-3\geq z>Y\geq 5.$ Then $z\geq21.$ Hence $Y\geq(z+3)/2\geq12.$ In view of $2\nmid Y, Y\geq 13.$ Then $\delta_{2}^{1/Y}\leq 3^{1/13}<3^{1/12},$ a contradiction.

This completes the proof of Theorem 1.3.

\section{Proof of Corollary 1.2 and Corollary 1.3}
$ $

We always assume that $(x, y, z)$ is a solution of the equation (1.1). In this section, We always assume that $r=2, n\geq 4.$
Recall that $$m=2^{\beta}j+e, n=2^{\alpha}i, 2\nmid ij, e\in\{1, -1\}, if 2\mid n.$$
$$m=2^{\alpha}i, n=2^{\beta}j+e, 2\nmid ij, e\in\{1, -1\}, if 2\nmid n.$$ Then $\alpha\geq1, \beta\geq 2.$

If $y>1.$ In view of the proof of Lemma 5.1, we may assume that $\alpha\geq 2.$ By [Mi09, Theorem 1.5], We only need to consider the case $2\alpha\neq \beta+1$ when $\alpha=2.$

\begin{lemma}\label{lemma1}
If $m\equiv 3(mod 4), n<64, 56n<m<(2n)^{10}$ and $(m, n)$ are not the case $$(\alpha, \beta)\in\{(3, 5), (4, 7), (5, 9)\} \quad \quad and \quad \quad j\equiv 1(mod 4).$$ Then conjecture 1.1 is true.
\end{lemma}

Proof. By Lemma 4.4, $y>1.$  In view of Lemma 2.2, 2.4, $2\mid x,  2\mid y.$
Since $(m, n)$ are not the case $$(\alpha, \beta)\in\{(3, 5), (4, 7), (5, 9)\} \quad \quad and \quad \quad j\equiv 1(mod 4).$$
Then $2\mid z$ by Lemma 2.5.
It follow from Lemma 3.2 that $$Y<\frac{\log 4c}{\log 4^{\alpha+1}}\leq\frac{\log 4c}{\log 64} \eqno{(7.1)}$$
In view of lemma 3.6(i), we may assume that $\Delta>0.$ Since $2\nmid XZ,$ then $\Delta=2(X-Z)\equiv 0(mod 4).$ By Lemma 6.2(ii), we have $Y\geq 5.$

If $m<(2n)^{2},$ then $Y<5$ by (7.1), a contradiction. Hence $m\geq(2n)^{2}>1.5n.$ By Lemma 6.1, $X<Z<Y.$ In view of Lemma 6.2(ii), $Y\geq 11.$
Since $$a^{2}+b^{2}=c^{2}, a^{2X}+b^{2Y}=c^{2Z},$$
Then $$\max\{a, b\}<c<\min\{a^{2}, b^{2}\}, \max\{a^{X}, b^{Y}\}<c^{Z}<\min\{a^{2X}, b^{2Y}\}.$$
Hence $c^{Z}<a^{2X}<c^{2X},$ then $Z<2X.$ It follow from Lemma 3.6 that
$$\frac{\log m}{\log n}<\Delta=2(Z-X)<Z<\frac{\log b}{\log c}Y+\frac{\log 2}{2\log c} \eqno{(7.2)}$$
If $m<(2n)^{4},$ then $Y<10,$ a contradiction. Hence $m>(2n)^{4}.$ Then $\Delta>4$ by (7.2), that is $\Delta\geq8$ by $4|\Delta.$ Then $Z\geq 9.$ By Lemma 6.2(ii), $Z\geq 11.$ Noting that $$f(m)=(\log(2mn))/(\log(m^{2}+n^{2}))=(\log b)/(\log c)$$ is a decreasing function of $m.$

If $m<(2n)^{6},$ then $Y<15.$ Hence $Y\leq 13.$ Then $$Z<\frac{5}{8}Y+\frac{\log 2}{2\log c}<9,$$ a contradiction.
Hence $m>(2n)^{6}.$

If $m<(2n)^{8},$ then $Y<19.$ That is $Y\leq17.$ Then $Z<\frac{7}{12}Y+\frac{\log 2}{2\log c}<10,$ a contradiction.
Hence $m>(2n)^{8}.$ Then $\Delta>8.$ So $\Delta\geq12,$ that is $Z\geq 13$ by (6.2).

If $m<(2n)^{10},$ then $Y\leq23.$ Hence $Z<\frac{9}{16}Y+\frac{\log 2}{2\log c}<13,$ a contradiction.

This completes the proof of Lemma 7.1. \hfill$\Box$\\

Let $q_{i}$ be the denominator of the $i-th$ convergent in the simple continued fraction expansion of $\gamma\in R,$
$a_{i}$ be the $i-th$ partial quotient of $\gamma\in R.$

\begin{lemma}\label{lemma1}~\cite{MiYW14}
Let $F_{n}, n\geq 0$ be the Fibonacci sequence given by $F_{0}=0, F_{1}=1, F_{n+2}=F_{n+1}+F_{n}, n\geq 0.$ We have $q_{k}\geq F_{k+1}.$
\end{lemma}

\begin{lemma}\label{lemma1}~\cite{Ka10}
Let $\gamma=\log a/\log c,$ then there exists a non-negative integer $s$ with $4\leq q_{s}<2521\log c$ such that
$$a_{s+1}+2>\frac{a^{q_{s}}\log c}{bq_{s}}.$$
\end{lemma}

\begin{lemma}\label{lemma1}
Let $\gamma=\log a/\log c,$ then $$2q_{s+1}>a^{q_{s}-2}.$$
\end{lemma}

Proof. By 7.2, there exists $s\in Z^{+}$ with $4\leq q_{s}<2521\log c$
such that $$a_{s+1}+2>\frac{a^{q_{s}}\log c}{bq_{s}}.$$
By $q_{s+1}=a_{s+1}q_{s}+q_{s-1},$
$$\frac{q_{s+1}-q_{s-1}}{q_{s}}+2>\frac{a^{q_{s}}\log c}{bq_{s}}.$$
That is $q_{s+1}-q_{s-1}+2q_{s}>\frac{a^{q_{s}}\log c}{b}.$
By $q_{s+1}\geq q_{s}\geq q_{s-1}, 2q_{s+1}>\frac{a^{q_{s}}\log c}{b}.$
Since $$a^{2}+b^{2}=c^{2},$$ then $$\max\{a, b\}<c<\min\{a^{2}, b^{2}\}.$$
Hence $2q_{s+1}>a^{q_{s}-2}$

This completes the proof of Lemma 7.4. \hfill$\Box$\\

\textbf{Proof of Corollary 1.2.} By Lemma 2.2, $2|x.$ If $y>1,$ in view of the proof of Lemma 5.1, we may assume that $\alpha\geq 2.$ In view of [Mi09, Theorem 1.5], We only need to consider the case $2\alpha\neq\beta+1$ when $\alpha=2.$ Since $m\equiv 3(mod 4)$ and $(m, n)$ not the case $$2\alpha=\beta+1, j\equiv 1(mod 4), \alpha\geq 3.$$ Then $j\neq e.$ By Lemma 2.5, $2\mid z.$

(1)\ We divide the argument into two cases:
Case 1: $m>56n.$ By Lemma 4.4, $y>1.$ In view of Lemma 2.4, $2|y.$
If $n\notin\{28, 44, 52, 56, 64\},$ the conjecture1.1 is true by Theorem 1.3. For $n=64,$ the conjecture is true by Theorem 1.2.
For $n\in\{28, 44, 52, 56\},$ by Lemma 7.1, we may assume that $m>(2n)^{10}.$ Then $$\log c>\log(m^{2})>\log((2n)^{20})>20\log n \eqno{(7.3)}$$
Then $$\nu=1+\frac{\log2}{\log(mn)}<1+\frac{\log2}{\log(4\times8^{10})}=\frac{33}{32}.$$
Hence $$n^{\nu}<4^{3}\leq 4^{\alpha+1}.$$

If $n=56,$ then $\alpha=3,$ then $$\eta=\frac{(22\nu+1)(\log 2)}{11\log(4^{\alpha+1}/n^{\nu})}<\frac{(22(33/32)+1)(\log 2)}{11\log(256/64)}<2.$$
If $n=28, 44,$ then $\alpha=2,$ hence $$\nu=1+\frac{\log2}{\log(mn)}<1+\frac{\log2}{\log(28\times56^{10})}<1.016.$$
$$\eta=\frac{(22\nu+1)(\log 2)}{11\log(4^{\alpha+1}/n^{\nu})}<\frac{(22\times1.016+1)(\log 2)}{11\log(64/47)}<5.$$
If $n=52,$ then $\alpha=2,$ hence $$\nu=1+\frac{\log2}{\log(mn)}<1+\frac{\log2}{\log(52\times104^{10})}<1.0138.$$
$$\eta=\frac{(22\nu+1)(\log 2)}{11\log(4^{\alpha+1}/n^{\nu})}<\frac{(22\times1.0138+1)(\log 2)}{11\log(64/54.915)}<9.6.$$
By (6.2), we have $\log c<\eta\log n<10\log n$ for $n\in\{28, 44, 52, 56\},$ which contradicts (7.3).

Case 2: \ $m<56n.$ Then $n\leq64, m<3584.$ If $y=1,$ then $3\leq z<x=q_{s}\leq 1534\log c\leq 25112.$ Hence $s\leq 24$ by Lemma 7.2. A simple computer calculation shows that for each of the pairs $(m, n)$ under consideration, the inequality $$a_{s+1}+2>\frac{a^{q_{s}}\log c}{bq_{s}}$$ dose not hold.

Hence $y>1,$ then $2|y$ by Lemma 2.4. In view of Lemma 3.2, $$Y<\frac{\log(4c)}{2(\alpha+1)\log 2}<4.27.$$ Hence $Y\leq 3$ by Lemma 2.6.
If $1.5n<m<56n,$ then $X<Z<Y$ by Lemma 6.1, that is $X<1<3,$ a contradiction.
If $m<1.5n,$ by Lemma 3.3, $$Y<\frac{\log(4c)}{2(\alpha+1)\log 2}<2.62.$$ Hence $Y=1$ by Lemma 2.6.
Since $$a^{2X}+b^{2Y}=c^{2Z},$$ Then $$a^{X}<b^{2Y}.$$ Hence $$X<\frac{2\log b}{\log a}Y\leq 4Y=4,$$ then $X\leq 3.$ Suppose $X<Z,$ by Lemma 3.6, $Z<Y=1,$ a contradiction. Hence $Z\leq X,$ then $Z\leq 3.$ If $Z=3,$ then $X=3,$ hence $a^{6}+b^{2}=c^{6}$ contradict $a^{2}+b^{2}=c^{2}.$ Hence $Z=1.$ Then $X=1.$ Hence Conjecture 1.1 is true.

(2)\ By Lemma 4.4, $y>1.$ In view of Lemma 2.4, $2\mid y.$ Since $n>64, m>56n,$ then $$\nu=1+\frac{\log 2}{\log(mn)}<1+\frac{\log 2}{\log(56\times64^{2})}<1.0562.$$
and $$n^{\nu-1}=e^{(\nu-1)\log n}=e^{\frac{\log 2}{\log(mn)}\log n}<e^{\frac{\log 2}{2}}<1.5.$$
By $n_{(2)}>\sqrt{n}, n=2^{\alpha}i,$ then $4^{\alpha+1}>4n.$
then $$\frac{4^{\alpha+1}}{n^{\nu}}>\frac{4n}{nn^{\nu-1}}=\frac{4}{n^{\nu-1}}>2.666,$$
Hence $$\eta=\frac{(22\nu+1)(\log 2)}{11\log(4^{\alpha+1}/n^{\nu})}<\frac{(22\times1.0562+1)(\log 2)}{11\log 2.666}<1.6.$$
By (6.2), we have $\log c<\eta\log n.$ Then $\log c<2\log n,$ a contradiction.

This completes the proof of Corollary 1.2. \hfill$\Box$\\

\textbf{Proof of Corollary 1.3.} We have $2\mid x$ by Lemma 2.2. Firstly we consider the case $y=1.$
If $m>56n$ or
\  \begin{displaymath}\begin{cases}
    \ u\geq94,  &  \text{$n=6^{u}$}, \\
    \ u\geq40,  &  \text{$n=10^{u}$}, \\
 \end{cases}                \end{displaymath}
Then $\max\{n_{(3)}, n_{(5)}\}\geq (1534\log n)^{1.5}t(n)\sqrt{n}.$ Hence the equation has no solution by Lemma 4.4, 4.7.

If $m<56n$ and
\  \begin{displaymath}\begin{cases}
    \ u<94,  &  \text{$n=6^{u}$}, \\
    \ u<40,  &  \text{$n=10^{u}$}, \\
 \end{cases}                \end{displaymath}
Then $3\leq z<x=q_{s}\leq 1534\log c\leq 529092.$ By Lemma 7.2, $s\leq 32.$ A simple computer calculation shows that for each of the pairs $(m, n)$ under consideration, the inequality $$a_{s+1}+2>\frac{a^{q_{s}}\log c}{bq_{s}}$$ dose not hold, hence the equation has no solution.

Hence $y>1.$ Then $2\mid y$ by Lemma 2.4. Noting that $$m=2^{\beta}j+e=2^{\beta}j-1, n=2^{\alpha}i, 2\nmid ij.$$

(1)\ $(m, n)$ are not the case $$2\alpha=\beta+1, j\equiv 1(mod 4), \alpha\geq 3.$$ Then $2\mid z$ by Lemma 2.5. Noting that $$2n_{(2)}n_{(3)}n_{(5)}\geq 2^{s(n)}\sqrt{n}.$$ Then Conjecture 1.1 is true by Theorem 1.3.

(2)\ $(m, n)$ are the case $$2\alpha=\beta+1, j\equiv 1(mod 4), \alpha\geq 3.$$ Then $n\geq 8.$ If $2\mid z,$ it seem to (1).
Then we only need to consider the case $2\nmid z.$ That is $2\nmid\Delta.$
If
\  \begin{displaymath}\begin{cases}
    \ u\geq 11,  &  \text{$n=6^{u}$}, \\
    \ u\geq 4,  &  \text{$n=10^{u}$}, \\
 \end{cases}                \end{displaymath}
Then $$2n_{(2)}n_{(3)}n_{(5)}\geq 2^{s(n)}\sqrt{n}, \quad \max\{n_{(3)}, n_{(5)}\}\geq 2\sqrt{n\log n}, \quad n_{(3)}n_{(5)}\geq\frac{3^{1/13}}{2^{\alpha+1}}n.$$
Hence Conjecture 1.1 is true by Theorem 1.3.

If
\  \begin{displaymath}\begin{cases}
    \ 3\leq u\leq 10,  &  \text{$n=6^{u}$}, \\
    \ u=3,  &  \text{$n=10^{u}$}, \\
 \end{cases}                \end{displaymath}
We divide the argument into two cases:

Case 1: $x>z.$ By Lemma 3.6(2), $$1<\frac{\log m}{\log n}<\Delta<\frac{2}{\log 3}\log F_{2}(K).$$
In view of (6.3), $$n^{2}_{(q)}\leq(m+1)\Delta.$$
Hence $$n^{2}_{(q)}\leq (m+1)\Delta<(m+1)\frac{2}{\log 3}\log F_{2}(K).$$
That is $$n^{2}_{(q)}<(m+1)\frac{2}{\log 3}\log F_{2}(K) \eqno{(7.3)}$$

Since $x>z,$ then $m<1.5n$ by Lemma 6.2. Noting that $(m, n)$ are the case $2\alpha=\beta+1, j\equiv 1(mod 4), \alpha\geq 3.$
Hence $$n=10^{3}, m=2^{5}j, j=33, 37, 41, 45, \quad or \quad n=6^{\alpha}, m=2^{2\alpha-1}j-1, j\equiv 1(mod 4), g\leq j\leq k,$$ where
$$\begin{tabular}{|r|r|r|r|r|r|r|r|r|}
\hline
$t$         & $3$    & $4$      & $5$      & $6$      & $7$      & $8$      & $9$       &$10$\\
\hline
$(g, k)$    &$(9, 9)$&$(13, 13)$&$(17, 21)$&$(25, 33)$&$(37, 49)$&$(53, 73)$&$(77, 113)$&$(117, 169)$\\
\hline
\end{tabular}$$

A simple calculation shows that for each of the pairs $(m, n)$ under consideration, (7.3) dose not hold.

Case 2: $x<z.$ Since $n_{(3)}n_{(5)}\geq\frac{3^{1/13}}{2^{\alpha+1}}n,$ by the proof of Theorem 1.3, we get a contradiction.

This completes the proof of Corollary 1.3. \hfill$\Box$\\

\section{Appendix: Estimates of $s(n), t(n)$ and one example}

In this section, we always assume that $r=2, n\geq 4.$ We will give estimates of $s(n), t(n)$ and an example of the case of $m=\prod_{i=1}^{k}p_{i}^{l_{i}}, p_{i}\equiv 1(mod 4).$
Noting that
\begin{eqnarray*}
s(n) &=&\frac{2\log(n^{11}-4n^{2})+5\log n+2\log(n^{11}/4)}{2(\log(n^{11}/4))\log((n^{11}-4n^{2})n^{2}/4)-2(\log2)\log n}\log n\\
&<&\frac{2\log(n^{11}-4n^{2})+5\log n+2\log(n^{11}/4)}{4(\log(n^{11}/4))-2(\log2)}\\
&<&\frac{2\log(n^{11})+5\log n+2\log(n^{11})-4\log 2}{4\log(n^{11})-10\log2}\\
&=&\frac{49\log n-4\log 2}{44\log n-10\log2}\\
&=&\frac{49}{44}+\frac{(157/22)\log2}{44\log n-10\log2}<1.2052\\
\end{eqnarray*}

Since $$1+\frac{2n^{2}}{6n+9}=\frac{2n^{2}+6n+9}{6n+9}<\frac{1}{2}n,$$ then
\begin{eqnarray*}
t(n) &=&\sqrt{\frac{1}{1534}+(1+\frac{4}{n})\log(1+\frac{2n^{2}}{6n+9})}\frac{\log(1+\frac{2n^{2}}{6n+9})}{(\log n)^{1.5}}\\
&<&\sqrt{\frac{1}{1534}+2\log(0.5n)}\frac{\log(0.5n)}{(\log n)^{1.5}}\\
&=&\sqrt{\frac{1}{1534}+2\log n-\log4}\frac{\log n-\log 2}{(\log n)^{1.5}}\\
&<&\sqrt{2\log n}\frac{\log n}{(\log n)^{1.5}}=\sqrt{2}\\
\end{eqnarray*}

\begin{lemma}\label{lemma1}~\cite{Mi11}

(1)\ Let $d\mid m+n.$ (i)If $d\equiv 7(mod 8),$ then $2\mid y.$

(ii)If $d\equiv 3(mod 8),$ then $2\mid z.$

(iii)If $d\equiv 5(mod 8)$ then $y\equiv z(mod 2).$

(2)\ Let $d\mid m-n.$ If $d\equiv \pm3(mod 8),$ then $y\equiv z(mod 2).$
\end{lemma}

\begin{exmp}\label{thm2}
If $(m, n)=(289\times 15233, 16)$ then Conjecture 1.1 is true.
\end{exmp}
If $(x, y, z)$ is a solution of the equation (1.1). Noting that $$m=17^{2}\times(17\times 16\times 56+1).$$
Then $m>56n,$ Hence $y>1.$ Since $17\mid gcd(m, n^{2}-1),$ then $2\mid x$ by Lemma 2.1.

In view of $m+n=3\times1467451,$ then $2\mid z$ by Lemma 8.1.

Since $m-n=7\times11\times57173,$ then $y\equiv z(mod 2)$ by Lemma 8.1, then $2\mid y.$

Since $m>56n,$ then $x<z$ by Lemma 6.1. It is similar to Theorem 1.2, Conjecture 1.1 is true.

\bigskip

{\bf Acknowledgement}\ \ Sincere thanks to Professor Pingzhi Yuan for his careful guidance.

This work was supported by the National Nature Science Foundatin of China, No.11671153, No.11971180.

 \end{document}